\documentclass[a4paper,12pt]{article}
\usepackage{amsbsy}
\usepackage{amsmath}
\usepackage{amsthm}
\usepackage{amstext}
\usepackage{amssymb}

\usepackage{xcolor}

\usepackage[T1]{fontenc}
\usepackage[latin1]{inputenc}
\usepackage{graphicx}

\PassOptionsToPackage{normalem}{ulem}
\usepackage{ulem}

\newtheorem{theorem}{Theorem}
\newtheorem{definition}{Definition}
\newtheorem{lemma}{Lemma}

\begin{document}

\author{Leonhard L{\"u}cken$^1$, Jan Philipp Pade$^2$, and Kolja Knauer$^3$}

\title{Classification of coupled dynamical systems with multiple delays: Finding the minimal number of delays}	


\maketitle

\noindent $^1$ {Weierstrass Institute, Mohrenstrasse 39, 10117 Berlin, Germany}\\
$^2$ {Humboldt University of Berlin, Institute of Mathematics, Unter den Linden 6, 10099 Berlin, Germany}\\
$^3$ {Universit{\'e} Montpellier 2, Laboratoire d'Informatique, de Robotique et de Micro{\'e}lectronique de Montpellier (LIRMM), Case Courrier 477, 161 rue Ada, 34095 Montpellier Cedex 5, France }

\begin{abstract}
In this article we study networks of coupled dynamical systems with
time-delayed connections. If two such networks hold different delays
on the connections it is in general possible that they exhibit different
dynamical behavior as well. We prove that for particular sets of
delays this is not the case. To this aim we introduce a componentwise
timeshift transformation (CTT) which allows to classify systems which
possess equivalent dynamics, though possibly different sets of connection
delays. In particular, we show for a large class of semiflows (including
the case of delay differential equations) that the stability of attractors
is invariant under this transformation. Moreover we show that each
equivalence class which is mediated by the CTT possesses a representative
system in which the number of different delays is not larger than
the cycle space dimension of the underlying graph. We conclude that
the 'true' dimension of the corresponding parameter space of delays
is in general smaller than it appears at first glance.
\end{abstract}

\section{Introduction}

Differential equations with time delays have been subject of intensive
research in the last decades. Recently, major impulses for theoretical investigations came 
from the challenge of understanding and modeling the human brain \cite{Izhikevich2008, Gosh2008, Deco2011}.
Here, time delays appear due to the finite speed of action potentials propagating through 
axons \cite{Manor1991, Wu2001, Campbell2007}. It was observed  that the 
qualitative behaviour of brain networks is closely related to their heterogeneous 
delay distributions. These were shown to play important roles in 
phenomena such as coherence in the resting state activity \cite{Deco2011,Gosh2008}, 
nonstationary bifurcations of equilibria \cite{Atay2004} and 
enhanced synchronizability \cite{Ernst1995,Dhamala2004}.
Further examples of delay-coupled systems are interacting lasers 
where delays appear due to the finite speed of light travelling 
through optical fibers \cite{Wunsche2005,Erneux2009,Soriano2013}. In population dynamics, 
delays correspond to maturing and gesturing times \cite{Kuang1993}, and 
in gene regulatory networks, delays represent 
the time the system needs to produce a protein \cite{Li2007,Danino2010}.
Thus, delays are often
included in order to account for the time a signal needs to propagate
from one node of the network to another or for the processing time
that it takes to emit a response to some input. An inevitable property
of real networks of interacting systems with delays is that, generically,
each delay time is different from any other. 

Nevertheless, there exist comparatively
few attempts to study systems with several different delays analytically.
An important reason for this is that their analysis is usually much
harder than for identical delays and techniques available 
for systems with one delay are in general not applicable in the case of multiple 
delays. In particular, the
analysis of spectra of solutions becomes more complicated since the 
characteristic quasi-polynomials involve several exponential terms \cite{Bellman1963, Hale1993}.
At least for coupled systems with two different delays 
some analytical results are available. For instance, 
Nussbaum \cite{Nussbaum1978} proved the existence of periodic
solutions in a system with two commensurable delays. 
Moreover some detailed studies on bifurcations in 
coupled systems with two different delay times were conducted 
\cite{Belair1994, Campbell1999, Shayer2000, Campbell2006, Guo2008}.
In particular, we mention Shayer et al. \cite{Shayer2000}, 
where the investigation initially assumes three different delays. 
In the course of the investigation, the authors discover 
that the stability of a steady state depends only on two 
values combined from these delays. This finding is a special case 
of the more general result that we present in this paper.
Higher order scalar systems with two delays are considered, 
e.g., by Gu et al. \cite{Gu2005}, who study stability crossing curves for 
this system. For the same system, Ruan and Wei \cite{Ruan2003} 
refine techniques for the determination of the roots of characteristic quasi-polynomials 
of single delay equations to treat the case of two delays.
Yanchuk and Giacomelli \cite{Yanchuk2014a} consider a scalar
system with two large delays and show that this system can
be described by a complex Ginzburg-Landau equation in the 
neighborhood of an equilibrium.

Finally, we want to mention studies of Hopfield neural networks with delayed
connections where many different delays are taken into account
\cite{Gopalsamy1994, Ye1995, Driessche1998, Lu2004}. Here,
the primary object of investigation are sufficient conditions for global 
convergence of the system towards a single steady state in order to obtain a 
well-defined method of input classification.

In this paper, the focus of our interest is a componentwise timeshift
transformation (CTT) which allows to change the interaction delays
while the dynamical properties of the system remain the same. Apart
of providing an intuition about the functionality of delays, the transformation
proves to be particularly useful in cases where it is possible to
achieve a smaller or in some sense preferable set of delays in the
transformed system.

This is especially apparent for the case when the nodes are coupled
in a single unidirectional ring. The corresponding special case of the CTT 
was utilized more or less explicitly already in \cite{Baldi1994, Mallet-Paret1996a, Perlikowski2010, Popovych2011}.
Therefore, let us illustrate the CTT for this example. In a general form, a ring of $N$ 
unidirectionally delay-coupled systems can be written as
\begin{equation}
\frac{\mathrm{d}}{\mathrm{dt}}x_{j}(t)=f_{j}\left(x_{j}(t),x_{j-1}(t-\tau_{j-1})\right),\ j=1,...,N,\label{eq:uniring-orig}
\end{equation}
where $j$ is considered modulo $N$ and $\tau_{j}\ge0$ are $N$
possibly different delays. Given a solution $x_{j}(t)$, $j=1,...,N$, we introduce
new variables
\begin{equation}
y_{j}\left(t\right):=x_{j}\left(t+\eta_{j}\right),\label{eq:new-variables-y}
\end{equation}
by shifting the nodes independently in time by certain amounts $\eta_{j}\ge0$. 
Formally, we find that the dynamics of the variables $y_{j}(t)$ obey 
\begin{equation}
\frac{\mathrm{d}}{\mathrm{dt}}y_{j}(t)=f_j\left(y_{j}(t),y_{j-1}(t-\tilde{\tau}_{j-1})\right),\ j=1,...,N,\label{eq:uniring-transf}
\end{equation}
with new delays $\tilde{\tau}_{j}=\tau_{j}+\eta_{j}-\eta_{j+1}$. 
In other words, finding timeshifts $\eta_j$ for given delays $\tilde{\tau}_j$ is equivalent to solving a system of linear equations.

\begin{figure}[t]%
\begin{centering}
\includegraphics[scale=0.25]{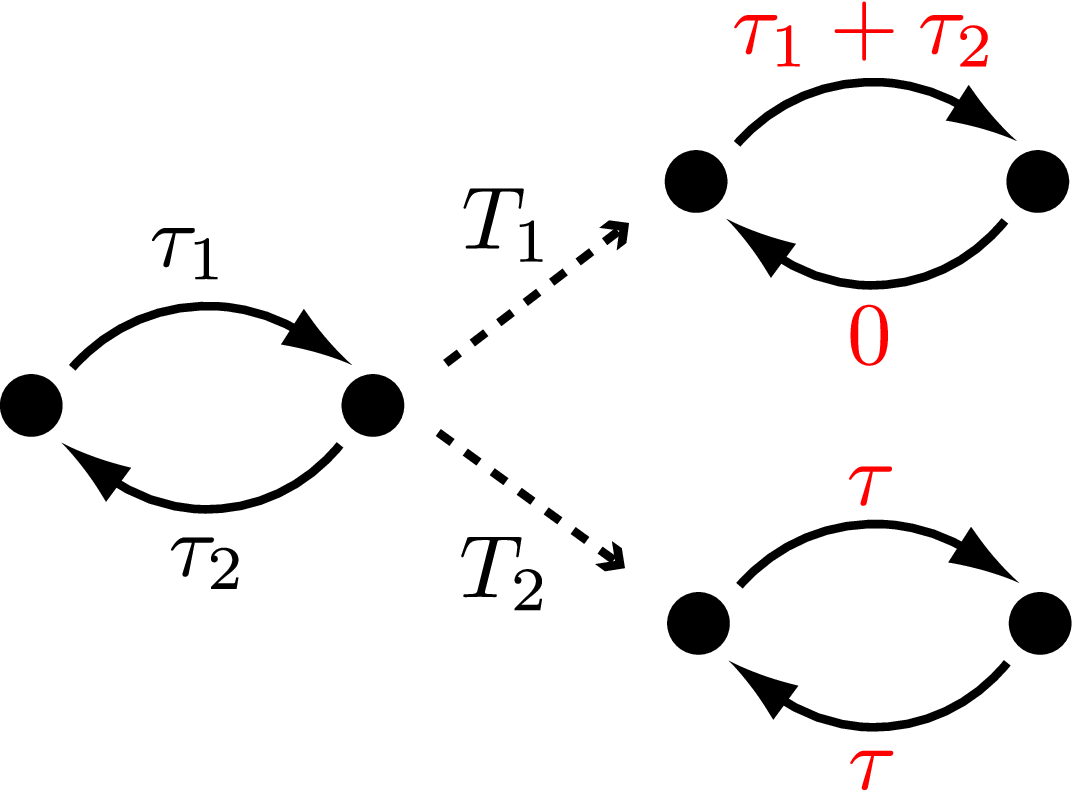}
\par\end{centering}

\caption{\label{fig:CTT-Example}Two simple examples of CT-Transformations
for a ring of two delay-coupled systems. The transformation $T_{1}$ reduces the number 
of delays from two to one and $T_{2}$ reveals a hidden symmetry in the 
system by adjusting the two delays to the same value $\tau=(\tau_{1}+\tau_{2})/2$.}
\end{figure}
\begin{figure}[t]%
\begin{centering}
\includegraphics[width=0.65\textwidth]{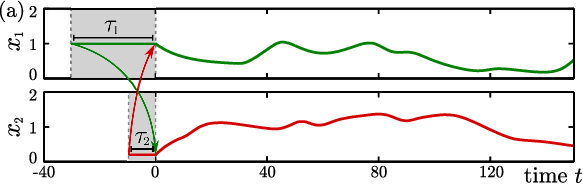}
\par\end{centering}

\begin{centering}
$\downarrow T$
\par\end{centering}

\begin{centering}
\includegraphics[width=0.65\textwidth]{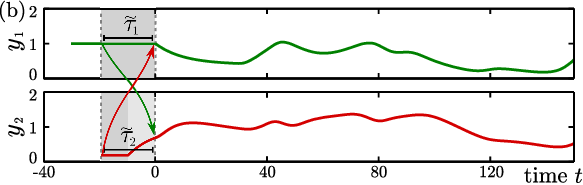}
\par\end{centering}

\begin{centering}
$\downarrow\tilde{T}$
\par\end{centering}

\centering{}\includegraphics[width=0.65\textwidth]{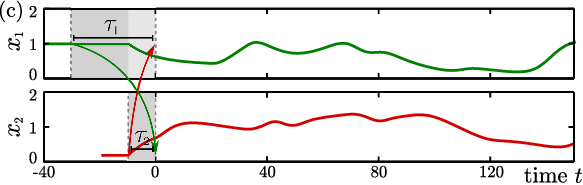}
\caption{\label{fig:CT-Transfo}CT-Transformation for a pair of bidirectionally
coupled Mackey-Glass systems \cite{Mackey1977}, cf. Fig.~\ref{fig:CTT-Example}. 
Plot (a) shows the original solution $\boldsymbol{x}(t)=(x_1(t),x_2(t))$ of (\ref{eq:1}).
Plot (b) shows the transformed solution $\boldsymbol{y}(t)=(y_1(t),y_2(t))$
of (\ref{eq:1-transformed}), which is shifted to the left by amounts
$\eta_{1}=0$ and $\eta_{2}=10$. Plot (c) shows the solution $\boldsymbol{x}(t+\eta_1)$
obtained after a reverse transformation of $\boldsymbol{y}(t)$ 
[see Sec.~\ref{sec:The-Componentwise-Timeshift} for details].}
\end{figure}

%
For the case of a ring with only two nodes, which is depicted in Fig.~\ref{fig:CTT-Example},
we illustrate the possible effect of the timeshift (\ref{eq:new-variables-y}) in Fig.~\ref{fig:CT-Transfo}.
The plots (a)--(c) show an initial piece of a solution and its change under the CTT
(two delay-coupled Mackey-Glass systems \cite{Mackey1977} were used to create this example).
The process of transformation is indicated by the symbols $T$ and $\tilde{T}$, which
are given a precise meaning in Sec.~\ref{sec:The-Componentwise-Timeshift}. 
Each transformation corresponds to a particular choice of the timeshifts in (\ref{eq:new-variables-y}).
Here, $T$ converts (\ref{eq:uniring-orig}) to (\ref{eq:uniring-transf})
and $\tilde{T}$ describes the reverse transformation from (\ref{eq:uniring-transf}) to (\ref{eq:uniring-orig}).
In (a) and (c), the arrows between the timetraces of both components indicate the delayed dependence
of $\dot{x}_1(t)$ on $x_2(t-\tau_2)$ and $\dot{x}_1(t)$ on $x_2(t-\tau_2)$,
similar in (b) for $y_1$ and $y_2$ for the transformed delays $\tilde{\tau}_{1}$ and $\tilde{\tau}_{2}$.
From (\ref{eq:new-variables-y}) it follows that the timetraces of the single components have exactly 
the same form in both systems. Both solutions only differ in the relative timeshifts between their components.
Baldi and Atiya \cite{Baldi1994} discovered that by choosing timeshifts
$\eta_{j} = \sum_{k=1}^{j-1} \tau_k$
one obtains the transformed delays
\begin{align}
\tilde{\tau}_{j} = \begin{cases}
0, &\text{for } 1\le j\le N-1,\\
\sum_{k=1}^{N}\tau_{k}, &\text{for } j=N.
\end{cases} \label{eq:Baldi-delays}
\end{align}
They used this simplified form to predict oscillatory 
behavior and bifurcations for a neural circuit exhibiting delayed excitatory and inhibitory
connections. Mallet-Paret \cite{Mallet-Paret1996a} utilized the same transformation to 
prove a Poincar{\'e}-Bendixson Theorem for monotone cyclic feedback systems with delays.
A slightly different set of timeshifts can be found in \cite{Perlikowski2010} and 
\cite{Popovych2011}, where 
the choice $\eta_{j} = (N-j)\bar{\tau} + \sum_{k=1}^{j} \tau_k$,  
with the mean delay $\bar{\tau}=\frac{1}{N}\sum_{k}\tau_{k}$,
was proposed. This leads to identical transformed delays
\begin{equation}
\tilde{\tau}_{j} \equiv \bar{\tau}. \label{eq:Perlikowski-delays}
\end{equation}
For a ring of identical nodes, i.e. $f_j \equiv f$ for all $j$,
it turns out that the obtained system is much more tractable
due to its rotational symmetry. One may say that this hitherto hidden 
symmetry was revealed by the CTT. See also Fig.~\ref{fig:CTT-Example} 
for an illustration of the two above mentioned timeshifts in a ring of two nodes.

For both variants, (\ref{eq:Baldi-delays}) and (\ref{eq:Perlikowski-delays}),
it is evident that if one is interested in the changes of dynamic behavior with respect
to the different delays $\tau_{j}$ it suffices to vary a single parameter,
i.e. the mean delay $\bar{\tau}=\frac{1}{N}\sum_{k}\tau_{k}$, instead
of the $N$ different parameters $\tau_{j}$, $j=1,...,N$. In other
words, the parameter space dimension is much smaller than it might
have appeared at first glance. Similar as for a ring, the CTT allows
to identify a canonical set of delay parameters in a general network
as we explain in Sec.~\ref{sec:Reduction-of-Delay-Parameters}.

Before doing so let us pose a question. What does the knowledge about
dynamical features of the transformed system (\ref{eq:uniring-transf})
really tell us about the dynamical features of the original system
(\ref{eq:uniring-orig})? It might seem intuitive that they are the
same since the timeshift $y_{j}\left(t\right)\mapsto y_{j}\left(t-\eta_{j}\right)$
inverses (\ref{eq:new-variables-y}). This is true for the particular
solution (\ref{eq:new-variables-y}) but, in general, the expression
$y_{j}\left(t-\eta_{j}\right)$ (with $\eta_{j}>0$) may not be defined
at all. This is because, in general, solutions of delay differential
equations (DDEs) cannot be continued backwards. Another subtlety arises
if one considers timeshifts $\eta_{j}$ leading to anticipating arguments,
that is negative delays. In this case (\ref{eq:uniring-transf}) has
fundamentally different properties from an ordinary DDE with positive
delays. For instance, initial value problems are ill-posed in general
\cite{Hale1993}. Even though we will restrict ourselves to the situation
where transformed delays are positive, a proper treatment of the above
question should introduce state spaces and flows to formulate and
compare the dynamical properties of the original and the transformed
system. Such a rigorous treatment of the CTT was not given in any 
of the above mentioned works \cite{Baldi1994, Mallet-Paret1996a, Perlikowski2010, Popovych2011}.
We will do this in a quite general setting using semidynamical
systems in Sec.~\ref{sec:The-Componentwise-Timeshift} where we
also give a rigorous definition of the CTT. In particular, the structural
similarity of the state spaces and the stability of invariant sets
is studied. In Sec.~\ref{sec:General-setup} we introduce notations
to describe general networks of delay coupled dynamical systems and
in Sec.~\ref{sub:Stability-of-Equilibria} we study the special
cases of equilibria and periodic orbits of DDEs where the equivalence
of the dynamics of the original and the transformed system is relatively
easy to show.

\section{Networks of delay coupled dynamical systems\label{sec:General-setup}}

To describe a general network of $N$ coupled systems we choose a
framework which enables us to account for multiple links between two
nodes holding different delays. That is, the coupling structure of
the network is assumed to be represented by a \emph{multidigraph}. This is
a set of node indices, ${\cal N}=\left\{ 1,...,N\right\} $, and a
set ${\cal E}$ of \emph{directed links}. Throughout the whole article
we assume that $({\cal N},{\cal E})$ is \emph{weakly connected}. This means that 
each node can be reached from any other node by traversing a sequence of
links, where each link may be traversed in arbitrary direction. For networks
with several connected components our results can be applied separately 
to each component. For each link $\ell\in{\cal E}$
the functions $s,t:{\cal E}\to{\cal N}$ assign its source $s(\ell)$ and its
target  $t(\ell)$. This means that the link $\ell$ connects the node $x_{s(\ell)}$
to the node $x_{t(\ell)}$. The delay time of $\ell$ is denoted by
$\tau(\ell)$. Note that there may indeed exist two links $\ell_{1}$
and $\ell_{2}$ in ${\cal E}$ with $s(\ell_{1})=s(\ell_{2})$ and
$t(\ell_{1})=t(\ell_{2})$ but $\tau(\ell_{1})\ne\tau(\ell_{2})$.
For node $x_{j}$, we introduce the set of its incoming links
as 
$$I_{j}=\left\{ \ell\in{\cal E}:\, t\left(\ell\right)=j\right\}.$$
Then, the dynamics of $x_{j}$ can be written as 
\begin{equation}
\frac{\mathrm{d}}{\mathrm{dt}}x_{j}(t)=f_{j}\left(\left(x_{s\left(\ell\right)}\left(t-\tau\left(\ell\right)\right)\right)_{\ell\in I_{j}}\right)\in\mathbb{R},\ j=1,...,N,\label{eq:1}
\end{equation}
where we assume $x_{j}(t)\in\mathbb{R}$ without loss of generality.
This notation allows to include a self-dependency of $x_{j}(t)$
via a link $\ell$ with $s(\ell)=t(\ell)=j$. It is also possible to 
have an instantaneous dependence by setting $\tau(\ell)=0$.
For (\ref{eq:1}), the introduction of new variables $y_{j}(t)$ as
in (\ref{eq:new-variables-y}) leads to the transformed system
\begin{equation}
\frac{\mathrm{d}}{\mathrm{dt}}y_{j}(t)=f_{j}\left(\left(y_{s\left(\ell\right)}\left(t-\tilde{\tau}\left(\ell\right)\right)\right)_{\ell\in I_{j}}\right),\ j=1,...,N,\label{eq:1-transformed}
\end{equation}
with modified delays
\begin{equation}
\tilde{\tau}\left(\ell\right)=\tau\left(\ell\right)-\eta_{t\left(\ell\right)}+\eta_{s\left(\ell\right)}.\label{eq:new-delays}
\end{equation}

\section{Spectrum of equilibria and periodic orbits\label{sub:Stability-of-Equilibria}}

An equilibrium point $\bar{\boldsymbol{x}}=\left(\bar{x}_{1},...,\bar{x}_{N}\right)\in\mathbb{R}^{N}$
of (\ref{eq:1}) is a point which satisfies
\begin{equation}
f_{j}\left(\left(\bar{x}_{s\left(\ell\right)}\right)_{\ell\in I_{j}}\right)=0,\ j=1,...,N. \label{eq:equilibrium-eqn}
\end{equation}
Obviously $\bar{\boldsymbol{x}}$ is also an equilibrium of (\ref{eq:1-transformed}).
In the following we assume that $f_{j}\in C^{1}.$ Then a characteristic
exponent $\lambda$ of $\bar{\boldsymbol{x}}$ in (\ref{eq:1}) corresponds to
an exponential solution $\boldsymbol{\xi}\left(t\right)=e^{\lambda t}\boldsymbol{\xi}_{0}\in\mathbb{C}^{N}$
of the variational equation 
\begin{equation}
\frac{\mathrm{d}}{\mathrm{dt}}{\xi}_{j}\left(t\right)=\sum_{\ell\in I_{j}}\partial_{s\left(\ell\right)}f_{j}\left(\left(\bar{x}_{s\left(\ell^{\prime}\right)}\right)_{\ell^{\prime}\in I_{j}}\right)\xi_{s(\ell)}\left(t-\tau\left(\ell\right)\right).  \label{eq:variational-eqn-fixed-point-original}
\end{equation}
If all of the characteristic exponents of $\bar{\boldsymbol{x}}$
possess negative real parts this assures that $\bar{\boldsymbol{x}}$
is stable, if at least one has positive real part then $\bar{\boldsymbol{x}}$
is unstable \cite{Hale1993}. In the transformed system (\ref{eq:1-transformed})
the timeshifted variation $\boldsymbol{\chi}\left(t\right)$, given
by $\chi_{j}\left(t\right)=\xi_{j}\left(t+\eta_{j}\right)$, is a
solution of the corresponding variational equation of $\bar{\boldsymbol{x}}$
\begin{equation}
\frac{\mathrm{d}}{\mathrm{dt}}\chi_{j}\left(t\right)=\sum_{\ell\in I_{j}}\partial_{s\left(\ell\right)}f_{j}\left(\left(\bar{x}_{s\left(\ell^{\prime}\right)}\right)_{\ell^{\prime}\in I_{j}}\right)\chi_{s(\ell)}\left(t-\tilde{\tau}\left(\ell\right)\right).  \label{eq:variational-eqn-fixed-point-transformed}
\end{equation}
Hence, the characteristic exponents of $\bar{\boldsymbol{x}}$ are
the same in (\ref{eq:1}) and (\ref{eq:1-transformed}).

Similarly, the stability of a periodic solution $\bar{\boldsymbol{x}}\left(t\right)=\bar{\boldsymbol{x}}\left(t+T\right)$
is determined by its Floquet exponents if their real parts are different
from zero. Each exponent $\lambda$ corresponds to a solution $\boldsymbol{\xi}\left(t\right)$
of the variational equation
\begin{equation}
\frac{\mathrm{d}}{\mathrm{dt}}\xi_{j}\left(t\right)=\sum_{\ell\in I_{j}}\partial_{s(\ell)}f_{j}\left(\left(\bar{x}_{s(\ell^\prime)}\left(t-\tau\left(\ell^{\prime}\right)\right)\right)_{\ell^{\prime}\in I_{j}}\right)\xi_{s(\ell)}\left(t-\tau\left(\ell\right)\right) \label{eq:variational-eqn-original}
\end{equation}
which has the form $\boldsymbol{\xi}\left(t\right)=e^{\lambda t}\boldsymbol{p}\left(t\right)$
with a periodic function $\boldsymbol{p}\left(t\right)=\boldsymbol{p}\left(t+T\right)$.
Again, the timeshifted solution $\boldsymbol{\chi}\left(t\right)$
with $\chi_{j}\left(t\right)=\xi_{j}\left(t+\eta_{j}\right)$ fulfills
the variational equation 
\begin{equation}
\frac{\mathrm{d}}{\mathrm{dt}}\chi_{j}\left(t\right)=\sum_{\ell\in I_{j}}\partial_{s(\ell)}f_{j}\left(\left(\bar{y}_{s(\ell^\prime)}\left(t-\tilde{\tau}\left(\ell^{\prime}\right)\right)\right)_{\ell^{\prime}\in I_{j}}\right)\chi_{s(\ell)}\left(t-\tilde{\tau}\left(\ell\right)\right)  \label{eq:variational-eqn-transformed}
\end{equation}
of the corresponding periodic solution $\bar{y}_{j}\left(t\right)=\bar{x}_{j}\left(t+\eta_{j}\right)$
in the transformed system. Let us summarize this section in a theorem:
\begin{theorem}
The following statements hold true:

\emph{(i)} A fixed point $\bar{\boldsymbol{x}}$ possesses the same characteristic
exponents in (\ref{eq:1}) and (\ref{eq:1-transformed}).

\emph{(ii)} A periodic solution $\bar{\boldsymbol{x}}\left(t\right)$ possesses
the same Floquet exponents in (\ref{eq:1}) as the corresponding transformed
solution $\bar{\boldsymbol{y}}\left(t\right)$ of (\ref{eq:1-transformed}).
\end{theorem}

\section{The componentwise timeshift transformation\label{sec:The-Componentwise-Timeshift}}

\subsection{Definitions}

In this section we introduce a rigorous formulation of the idea which
underlies the change of variables (\ref{eq:new-variables-y}). We
define the CTT in terms of the underlying infinite dimensional phase
spaces of (\ref{eq:1}) and (\ref{eq:1-transformed}).

First recall that a \emph{semidynamical system} (or a \emph{semiflow})
is a mapping 
\begin{eqnarray*}
\Phi:\left[0,\infty\right)\times X & \rightarrow & X,\\
\left(t,\boldsymbol{x}\right) & \mapsto & \Phi_{t}\left(\boldsymbol{x}\right),
\end{eqnarray*}
on a Banach space $X$ which fulfills:
\[
\begin{array}{ll}
\text{(i)} & \Phi_{0}=Id\\
\text{(ii)} & \Phi_{t+s}=\Phi_{t}\circ\Phi_{s},\,\forall t,s\geq0\\
\text{(iii)} & \Phi_{t}:X\rightarrow X\text{ is continuous for all }t\geq0.
\end{array}
\]

The state spaces for DDEs contain segments of functions,
which represent the history of the solution curve $\boldsymbol{x}\left(t\right)$.
For instance, in Fig.~\ref{fig:CT-Transfo} the shaded part of the 
timetrace in (a) corresponds to the initial segment of the depicted
solution.
For the original system (\ref{eq:1}), we choose the state space to be
\[
{\cal C}=\prod_{j=1}^{N}C\left(\left[-r_{j},0\right];\,\mathbb{R}\right),\ r_{j}=\max_{\ell\in O_{j}}\tau\left(\ell\right),
\]
where $O_{j}=\left\{ \ell\in{\cal E}:\, s\left(\ell\right)=j\right\} $
is the set of all outgoing links from the $j$-th node. Hence, the
value $r_{j}$ is the largest delay time on outgoing links of the
node $j$. This definition ensures, that the history for the $j$-th
component is available for all delayed arguments appearing on the
right hand side of (\ref{eq:1}). Similarly, we choose 
\[
\tilde{{\cal C}}=\prod_{j=1}^{N}C\left(\left[-\tilde{r}_{j},0\right];\,\mathbb{R}\right),\ \tilde{r}_{j}=\max_{\ell\in O_{j}}\tilde{\tau}\left(\ell\right)
\]
as state space for the transformed system (\ref{eq:1-transformed}).
Assuming that solutions exist for all future times (e.g. if $f_{j}$
are Lipschitz continuous), there exist semiflows 
\begin{eqnarray*}
\Phi:\left[0,\infty\right)\times{\cal C}\to{\cal C}, &  & \ \ \left(t,\boldsymbol{x}\right)\mapsto\Phi_{t}\left(\boldsymbol{x}\right),\\
\Psi:\left[0,\infty\right)\times\tilde{{\cal C}}\to\tilde{{\cal C}}, &  & \ \ \left(t,\boldsymbol{y}\right)\mapsto\Psi_{t}\left(\boldsymbol{y}\right),
\end{eqnarray*}
for (\ref{eq:1}) and (\ref{eq:1-transformed}). Now let us formulate
(\ref{eq:new-variables-y}) as a state space transformation $T:{\cal C}\to\tilde{{\cal C}}$.
We define $T$ componentwise for $j=1,...,N$, and pointwise for $\boldsymbol{x}_{0}\in{\cal C}$
and $t\in[-\tilde{r}_{j},0]$, as
\begin{equation}
T_{j}\left[\boldsymbol{x}_{0}\right]\left(t\right)=\begin{cases}
\left[\boldsymbol{x}_{0}\left(t+\eta_{j}\right)\right]_{j}, & \text{for } t\in\left[-\tilde{r}_{j},-\eta_{j}\right],\\
\left[\Phi_{t+\eta_{j}}\left(\boldsymbol{x}_{0}\right)\right]_{j}\left(0\right), & \text{for } t\in\left[-\min\left\{ \eta_{j},\tilde{r}_{j}\right\} ,0\right].
\end{cases}\label{eq:def-T}
\end{equation}
For illustration see Fig.~\ref{fig:CT-Transfo}(b), where the shading marks the 
segment $\boldsymbol{y}_0=T[\boldsymbol{x}_0]$ for the initial 
segment $\boldsymbol{x}_0$ indicated in (a). The lighter shading indicates 
the part defined by the second case of (\ref{eq:def-T}).
Let us show that (\ref{eq:def-T}) is well-defined. For this we need
to assure that $\boldsymbol{y}_{0}=T\left[\boldsymbol{x}_{0}\right]\in\tilde{{\cal C}}$.
This only requires that the term $[\boldsymbol{x}_{0}(t-\eta_{j})]_{j}$
appearing in (\ref{eq:def-T}) is defined for all $t+\eta_{j}$ with
$t\in[-\tilde{r}_{j},-\eta_{j}]$. For the case in which $[-\tilde{r}_{j},-\eta_{j}]$
is non-empty, this is equivalent to $t+\eta_{j}\ge-r_{j}$ for all
$t\in[-\tilde{r}_{j},-\eta_{j}]$. In order to show this, choose a
link $\ell\in O_{j}$ with maximal delay $\tilde{\tau}(\ell)=\tilde{r}_{j}$.
Then,
\begin{equation*}
r_{j}-\tilde{r}_{j} \ge \tau\left(\ell\right)-\tilde{\tau}\left(\ell\right)
 = \tau\left(\ell^{\prime}\right)-\left(\tau\left(\ell^{\prime}\right)-\eta_{t\left(\ell^{\prime}\right)}+\eta_{j}\right)
 \ge-\eta_{j}.
\end{equation*}
Therefore, $t+\eta_{j}\ge-\tilde{r}_{j}+\eta_{j}\ge-r_{j}$ for $t\in[-\tilde{r}_{j},-\eta_{j}]$
and (\ref{eq:def-T}) is well-defined.

Furthermore, one easily checks that 
\begin{equation}
T\circ\Phi_{t}=\Psi_{t}\circ T,\ \text{for all }t\geq0.\label{eq:commutation-T-flow}
\end{equation}
That is, $T$ transforms solutions of (\ref{eq:1}) into solutions of (\ref{eq:1-transformed}).
As inherited from the semiflow $\Phi$,
the transformation $T$ is neither injective nor surjective in general.
Therefore, one cannot expect a dynamical equivalence of (\ref{eq:1})
and (\ref{eq:1-transformed}) in the strict form of topological
conjugacy, i.e. $\Psi=h\circ\Phi\circ h^{-1}$ for some homeomorphism
$h$. Moreover, since $T$ is not surjective, (\ref{eq:commutation-T-flow}) 
does not even signify that $\Psi$ is properly semiconjugate to $\Phi$. However,
this "weak semiconjugacy" is mutual as we show in the following,
and this fact implies a strong equivalence as well.

\subsection{CT-equivalence}

Let us find a reverse transformation from (\ref{eq:1-transformed})
to (\ref{eq:1}). It should transform $\tilde{\tau}\left(\ell\right)$
back to $\tau\left(\ell\right)$ which leads to the natural definition
of reverse timeshifts 
\[
\tilde{\eta}_{j}=\bar{\eta}-\eta_{j}\ge0,\ \text{with }\bar{\eta}=\max_{1\le j\le N}\eta_{j}.
\]
Then, the reverse transformation $\tilde{T}:\tilde{{\cal C}}\to{\cal C}$
is given as
\begin{equation}
\tilde{T}_{j}\left[\boldsymbol{y}_{0}\right]\left(t\right)=\begin{cases}
\left(\Psi_{t+\tilde{\eta}_{j}}\left(\boldsymbol{y}_{0}\right)\right)_{j}\left(0\right), & t\in\left[-\min\left\{ \tilde{\eta}_{j},r_{j}\right\} ,0\right],\\
\left(\boldsymbol{y}_{0}\left(t+\tilde{\eta}_{j}\right)\right)_{j}\left(0\right), & t\in\left[-r_{j},-\tilde{\eta}_{j}\right].
\end{cases}\label{eq:def-T-tilde}
\end{equation}
Analogously to (\ref{eq:commutation-T-flow}), we have 
\begin{equation}
\tilde{T}\circ\Psi=\Phi\circ\tilde{T},\label{eq:commutation-T-tilde-flow}
\end{equation}
 and additionally $T$ and $\tilde{T}$ are reverse in the sense that
\begin{equation}
\tilde{T}\circ T=\Phi_{\bar{\eta}} \ \text{ and } \ T\circ\tilde{T}=\Psi_{\bar{\eta}}.\label{eq:reverseness-T-and-T-tilde}
\end{equation}
See, for example, Fig.~\ref{fig:CT-Transfo}(c), where the shaded region
indicates the initial segment $(\tilde{T}\circ T)[\boldsymbol{x}_0]$ of the solution 
$(\tilde{T}\circ T)[\boldsymbol{x}](t) = \boldsymbol{x}(t+\bar{\eta})$, with $\bar{\eta}=10$.
The lighter shaded region indicates where the second case of (\ref{eq:def-T-tilde}) 
takes effect.

The equations (\ref{eq:commutation-T-flow}), (\ref{eq:commutation-T-tilde-flow})
and (\ref{eq:reverseness-T-and-T-tilde}) constitute a notion of equivalence
which can be applied not only to semidynamical systems which stem
from systems of DDEs, but as well to systems such as delay coupled
PDEs or more generally to abstract evolution equations. Therefore
we define in a more abstract fashion:
\begin{definition}
\label{def:CT-equivalence} Two semidynamical systems $\Phi:\left[0,\infty\right)\times X\to X$
and $\Psi:\left[0,\infty\right)\times Y\to Y$ are called CT-equivalent
if there exist mappings $T:X\to Y$, $\tilde{T}:Y\to X$ and $\bar{\eta}\ge0$
such that
\begin{align*}
T\circ\tilde{T}=\Psi_{\bar{\eta}}, & \ T\circ\Phi_{t}=\Psi_{t}\circ T,\\
\tilde{T}\circ T=\Phi_{\bar{\eta}}, & \ \tilde{T}\circ\Psi_{t}=\Phi_{t}\circ\tilde{T}.
\end{align*}

\end{definition}
In the following, the analysis is carried out in the general setting
on Banach spaces $X$ and $Y$ but we keep in mind that the results
apply to the case of DDEs on the spaces $X=\mathcal{C}$ and $Y=\tilde{\mathcal{C}}$
as introduced above. We still use the term CT-equivalence for the general formulation though,
since most probably the equivalent systems $\Phi$ and $\Psi$ will
in practice be connected by a transformation, which resembles (\ref{eq:def-T}).

\subsection{Dynamical invariants of CT-equivalent semidynamical systems}

In this section, we derive some properties of CT-equivalent systems
in the sense of Definition~\ref{def:CT-equivalence}. We show that
their state spaces hold a structure of corresponding strongly invariant
sets. Then, assuming that the CTTs $T$ and $\tilde{T}$ are
Lipschitz continuous, we prove that stability properties are preserved
as well. For the DDE systems (\ref{eq:1}) and (\ref{eq:1-transformed})
this is the case if all $f_{j}$ are Lipschitz continuous. Throughout,
we assume that $X$ and $Y$ are Banach spaces. Before stating the
main results of this section (Theorems~\ref{thm:correspondence}
and \ref{thm:stability}), we give some definitions:

A set $A\subseteq X$ is called \emph{positively invariant} under
$\Phi$ if for any $t\geq0$:
\[
\Phi_{t}\left(A\right)\subseteq A.
\]
$A\subseteq X$ is called \emph{invariant }under $\Phi$ if for any
$t\ge0$:
\[
\Phi_{t}\left(A\right)=A.
\]
For any set $A\subseteq X$, we define its \emph{strongly invariant
hull} $H_{\Phi}\left(A\right)$ as the set
\[
H_{\Phi}\left(A\right):=\left\{ \boldsymbol{x}\in X\,\mid\,\exists t_{1},t_{2}\ge0,\,\hat{\boldsymbol{x}}\in A:\,\Phi_{t_{1}}\left(\hat{\boldsymbol{x}}\right)=\Phi_{t_{2}}\left(\boldsymbol{x}\right)\right\} .
\]
A \emph{strongly invariant set}\textbf{\emph{ }}$A\subseteq X$ is
a set that coincides with its strongly invariant hull, i.e.
\[
H_{\Phi}\left(A\right)=A.
\]
The class of strongly invariant sets is denoted by
\[
\mathrm{sis}\left(\Phi\right)=\left\{ A\subseteq X\,\mid\, A=H_{\Phi}\left(A\right)\right\} .
\]
Note that positive invariance is implied by both, invariance and strong
invariance. But between invariance and strong invariance there holds
no implication. Of course a strongly invariant set always contains
a maximal invariant set which might be empty.
\begin{theorem}
\label{thm:correspondence} Let $\Phi:\left[0,\infty\right)\times X\to X$
and $\Psi:\left[0,\infty\right)\times Y\to Y$ be CT-equivalent. Then,

\emph{(i)} for each (positively) $\Phi$-invariant set $A\in X$, the set
$T\left[A\right]\in Y$ is (positively) $\Psi$-invariant,

\emph{(ii)} there is a one-to-one correspondence between (strongly) invariant
sets of $\Phi$ and $\Psi$. \end{theorem}
\begin{proof}
Ad~(i): Let $A$ be positively invariant, $\boldsymbol{x}\in A$
and $\boldsymbol{y}=T\left[\boldsymbol{x}\right]$. Then, $\Psi_{t}\left(\boldsymbol{y}\right)=T\left[\Phi_{t}\left(\boldsymbol{x}\right)\right]\in T\left[A\right]$.
Hence $T\left[A\right]$ is positively invariant. \\
If $A$ is invariant, then for each $\boldsymbol{x}\in A$ and each
$t\in\left[0,\infty\right)$ there is an $\boldsymbol{x}_{-t}\in A$
such that $\Phi_{t}\left(\boldsymbol{x}_{-t}\right)=\boldsymbol{x}$
and correspondingly for each $\boldsymbol{y}=T\left[\boldsymbol{x}\right]\in T\left[A\right]$
there is $\boldsymbol{y}_{-t}=T\left[\boldsymbol{x}_{-t}\right]$
with 
\[
\Psi_{t}\left(\boldsymbol{y}_{-t}\right)=T\left[\Phi_{t}\left(\boldsymbol{x}_{-t}\right)\right]=T\left[\boldsymbol{x}\right]=\boldsymbol{y}.
\]

Ad~(ii): Let $A\in\mathrm{sis}\left(\Phi\right)$. We define a corresponding
set $\tilde{A}=H_{\Psi}\left(T\left[A\right]\right)\in\mathrm{sis}\left(\Psi\right)$.
Correspondence for strongly invariant sets is proven via
\[
(H_{\Phi}\circ\tilde{T})\circ\left(H_{\Psi}\circ T\right)=\mathrm{id}_{\mathrm{sis}\left(\Phi\right)}\text{ and }\left(H_{\Psi}\circ T\right)\circ(H_{\Phi}\circ\tilde{T})=\mathrm{id}_{\mathrm{sis}\left(\Psi\right)},
\]
where, by symmetry, it suffices to show only one equality. Let $\boldsymbol{y}\in\tilde{A}$,
then there exist $t_{1},t_{2}\ge0$ and $\hat{\boldsymbol{y}}=T\left[\hat{\boldsymbol{x}}\right]\in T\left[A\right]$
such that $\Psi_{t_{2}}\left(\boldsymbol{y}\right)=\Psi_{t_{1}}\left(\hat{\boldsymbol{y}}\right)\in\tilde{A}$.
Thus,
\begin{eqnarray*}
\Phi_{t_{2}}\left(\tilde{T}\left[\boldsymbol{y}\right]\right) & = & \tilde{T}\left[\Psi_{t_{2}}\left(\boldsymbol{y}\right)\right]=\tilde{T}\left[\Psi_{t_{1}}\left(\hat{\boldsymbol{y}}\right)\right]\\
 & = & \tilde{T}\left[\Psi_{t_{1}}\left(T\left[\hat{\boldsymbol{x}}\right]\right)\right]=\tilde{T}\left[T\left[\Phi_{t_{1}}\left(\hat{\boldsymbol{x}}\right)\right]\right]\\
 & = & \Phi_{t_{1}+\bar{\eta}}\left(\hat{\boldsymbol{x}}\right)\in H_{\Phi}\left(A\right)=A.
\end{eqnarray*}
Hence, $\tilde{T}\left[\boldsymbol{y}\right]\in A$. That is, $\tilde{T}\left[\tilde{A}\right]\subseteq A$
and $H_{\Phi}\left(\tilde{T}\left[\tilde{A}\right]\right)\subseteq A$.
We have also shown that for $\hat{\boldsymbol{x}}\in A$, $\Phi_{t_{1}+\bar{\eta}}\left(\hat{\boldsymbol{x}}\right)\in\tilde{T}\left[\tilde{A}\right]$
and therefore $\hat{\boldsymbol{x}}\in H_{\Phi}\left(\tilde{T}\left[\tilde{A}\right]\right)$.
This yields $A\subseteq H_{\Phi}\left(\tilde{T}\left[\tilde{A}\right]\right)\subseteq A$,
i.e. $A=H_{\Phi}\left(\tilde{T}\left[\tilde{A}\right]\right)$ and
$\left(H_{\Phi}\circ\tilde{T}\right)\circ\left(H_{\Psi}\circ T\right)=\mathrm{id}_{\mathrm{sis}\left(\Phi\right)}$.

For invariant sets, the one-to-one correspondence is mediated directly
via $T$ and $\tilde{T}$.\end{proof}
\begin{definition}
The maximal Lyapunov exponent (MLE) of a point $\boldsymbol{x}\in X$
with respect to a semidynamical system $\Phi:\left[0,\infty\right)\times X\to X$
on a Banach space $X$ is defined as
\[
\lambda\left(\boldsymbol{x}\right):=\limsup_{t\to\infty}\limsup_{\left|\xi\right|\searrow0}\frac{1}{t}\ln\left(\frac{\left|\Phi_{t}\left(\boldsymbol{x}+\xi\right)-\Phi_{t}\left(\boldsymbol{x}\right)\right|}{\left|\xi\right|}\right)\in\left[-\infty,\infty\right].
\]
The MLE of $\boldsymbol{x}\in A\subseteq X$ with respect to the set
$A$ is defined as 
\begin{eqnarray*}
\lambda\left(x,A\right) & := & \limsup_{t\to\infty}\limsup_{\left|\xi\right|\searrow0}\min_{a\in A}\frac{1}{t}\ln\left(\frac{\left|\Phi_{t}\left(\boldsymbol{x}+\xi\right)-a\right|}{\left|\xi\right|}\right)\\
 & = & \limsup_{t\to\infty}\limsup_{\left|\xi\right|\searrow0}\frac{1}{t}\ln\left(\frac{\mathrm{dist}\left(\Phi_{t}\left(\boldsymbol{x}+\xi\right),A\right)}{\left|\xi\right|}\right).
\end{eqnarray*}
The MLE of a set $A\in X$ is defined as $\lambda\left(A\right)=\sup_{\boldsymbol{x}\in A}\lambda\left(\boldsymbol{x},A\right).$\end{definition}
\begin{theorem}
\label{thm:stability} Let the CT-transformations $T$ and $\tilde{T}$
be Lipschitz-continuous and let $X$, $Y$ be Banach spaces. Then,

\emph{(i)} corresponding (strongly) invariant sets of $\Phi$ and $\Psi$
possess the same maximal Lyapunov exponents (MLEs),

\emph{(ii)} for each positively invariant set $A\in X$, the set $T\left[A\right]\in Y$
has the same type of stability.
\end{theorem}
Claim (i) is a direct consequence of the following Lemma.
\begin{lemma}
\label{lem:Lyap-ineq} Let $X$, $Y$ be Banach spaces and let $T$
and $\tilde{T}$ be Lipschitz-continuous with constants $L_{T}$ and
$L_{\tilde{T}}$. Then, we have $\lambda\left(\boldsymbol{x}\right)\le\lambda\left(T\left[\boldsymbol{x}\right]\right)\le\lambda\left(\Phi_{\bar{\eta}}\left(\boldsymbol{x}\right)\right)$
for all $\boldsymbol{x}\in X$.\end{lemma}
\begin{proof}
For each $\boldsymbol{x},\boldsymbol{\chi}\in{\cal C}$, $t\ge\bar{\eta}$:
\begin{eqnarray*}
 &  & \frac{1}{t}\ln\left(\frac{\left|\Phi_{t}\left(\boldsymbol{x}\right)-\Phi_{t}\left(\boldsymbol{x}+\boldsymbol{\chi}\right)\right|}{\left|\boldsymbol{\chi}\right|}\right)\\
 & = & \frac{1}{t}\ln\left(\frac{\left|\Phi_{t-\bar{\eta}}\left(\tilde{T}\circ T\left[\boldsymbol{x}\right]\right)-\Phi_{t-\bar{\eta}}\left(\tilde{T}\circ T\left[\boldsymbol{x}+\boldsymbol{\chi}\right]\right)\right|}{\left|\boldsymbol{\chi}\right|}\right)\\
 & = & \frac{1}{t}\ln\left(\frac{\left|\tilde{T}\circ\Psi_{t-\bar{\eta}}\left(T\left[\boldsymbol{x}\right]\right)-\tilde{T}\circ\Psi_{t-\bar{\eta}}\left(T\left[\boldsymbol{x}\right]+\boldsymbol{\xi}\right)\right|}{\left|\boldsymbol{\chi}\right|}\right),
\end{eqnarray*}
with $\left|\boldsymbol{\xi}\right|=\left|T\left[\boldsymbol{x}+\boldsymbol{\chi}\right]-T\left[\boldsymbol{x}\right]\right|\le L_{T}\left|\boldsymbol{\chi}\right|$,
\begin{eqnarray*}
... & \le & \frac{1}{t}\ln\left(\frac{L_{\tilde{T}}\left|\Psi_{t-\bar{\eta}}\left(T\left[\boldsymbol{x}\right]\right)-\Psi_{t-\bar{\eta}}\left(T\left[\boldsymbol{x}\right]+\boldsymbol{\xi}\right)\right|}{\left(\left|\boldsymbol{\xi}\right|/L_{T}\right)}\right)\\
 & = & \frac{1}{t}\ln\left(\frac{\left|\Psi_{t-\bar{\eta}}\left(T\left[\boldsymbol{x}\right]\right)-\Psi_{t-\bar{\eta}}\left(T\left[\boldsymbol{x}\right]+\boldsymbol{\xi}\right)\right|}{\left|\boldsymbol{\xi}\right|}\right)+\frac{1}{t}\ln\left(L_{T}L_{\tilde{T}}\right)
\end{eqnarray*}
Therefore, 
\[
\frac{1}{t}\ln\left(\frac{\left|\Phi_{t}\left(\boldsymbol{x}\right)-\Phi_{t}\left(\boldsymbol{x}+\boldsymbol{\chi}\right)\right|}{\left|\boldsymbol{\chi}\right|}\right)\le\frac{1}{t}\ln\left(\frac{\left|\Psi_{t-\bar{\eta}}\left(T\left[\boldsymbol{x}\right]\right)-\Psi_{t-\bar{\eta}}\left(T\left[\boldsymbol{x}\right]+\boldsymbol{\xi}\right)\right|}{\left|\boldsymbol{\xi}\right|}\right),
\]
and, thus, 
\[
\lambda\left(\boldsymbol{x}\right)\le\lambda\left(T\left[\boldsymbol{x}\right]\right).
\]
The same reasoning for $\boldsymbol{y}=T\left[\boldsymbol{x}\right]$
gives
\[
\lambda\left(\boldsymbol{x}\right)\le\lambda\left(T\left[\boldsymbol{x}\right]\right)\le\lambda\left(\Phi_{\bar{\eta}}\left(\boldsymbol{x}\right)\right).
\]

\end{proof}

\begin{proof}
(of Theorem~\ref{thm:stability}) Lemma~\ref{lem:Lyap-ineq} implies
claim~(i). To see that consider two corresponding sets $A\in\mathrm{sis}\left(\Phi\right)$
and $\tilde{A}\in\mathrm{sis}\left(\Psi\right)$. Then it is impossible
that $\lambda\left(A\right)>\lambda(\tilde{A})$ since
for any $\boldsymbol{x}\in A$ there exists $\boldsymbol{y}=T\left[\boldsymbol{x}\right]\in\tilde{A}$
such that $\lambda\left(\boldsymbol{x}\right)\le\lambda\left(\boldsymbol{y}\right)$.
The same holds vice versa and therefore $\lambda\left(A\right)=\lambda(\tilde{A})$.
Similarly, one shows this for the case of invariant sets.

Ad~(ii): Let $A\in{\cal C}$ be an (asymptotically) stable positively
invariant set. We show that if $A$ is (asymptotically) stable so
is $T\left[A\right]$. Let $\varepsilon>0$, $\boldsymbol{y}_{0}=T\left[\boldsymbol{x}_{0}\right]\in T\left[A\right]$
and $\boldsymbol{\xi}_{0}\in\tilde{{\cal C}}$ a small initial perturbation,
i.e. $\left|\boldsymbol{\xi}_{0}\right|\le\delta$ with $\delta=\delta\left(\varepsilon\right)>0$
to be specified later. Define the perturbed solution
\[
\tilde{\boldsymbol{y}}_{t}=\Psi_{t}\left(\boldsymbol{y}_{0}+\boldsymbol{\xi}_{0}\right)=\boldsymbol{y}_{t}+\boldsymbol{\xi}_{t},
\]
where $\boldsymbol{y}_{t}=\Psi_{t}\left(\boldsymbol{y}_{0}\right)$
is the unperturbed solution. Consider
\begin{eqnarray*}
\tilde{\boldsymbol{x}}_{t}: & = & \tilde{T}\left[\tilde{\boldsymbol{y}}_{t}\right]=\tilde{T}\left[\Psi_{t}\left(\boldsymbol{y}_{0}+\boldsymbol{\xi}_{0}\right)\right]\\
 & = & \Phi_{t}\left(\tilde{T}\left[\boldsymbol{y}_{0}+\boldsymbol{\xi}_{0}\right]\right)\\
 & = & \Phi_{t}\left(\tilde{T}\circ T\left[\boldsymbol{x}_{0}\right]+\left(\tilde{T}\left[\boldsymbol{y}_{0}+\boldsymbol{\xi}_{0}\right]-\tilde{T}\left[\boldsymbol{y}_{0}\right]\right)\right)\\
 & = & \Phi_{t}\left(\boldsymbol{x}_{\bar{\eta}}+\boldsymbol{\chi}_{\bar{\eta}}\right),
\end{eqnarray*}
with 
\begin{eqnarray*}
\left|\boldsymbol{\chi}_{\bar{\eta}}\right| & = & \left|T\left[\boldsymbol{y}_{0}+\boldsymbol{\xi}_{0}\right]-\tilde{T}\left[\boldsymbol{y}_{0}\right]\right|\\
 & \le & L_{\tilde{T}}\cdot\delta,
\end{eqnarray*}
where $L_{\tilde{T}}$ is the Lipschitz constant of $\tilde{T}$.
Since $\boldsymbol{x}_{\bar{\eta}}\in A$ and $A$ is (asymptotically)
stable, we can find a $\delta=\delta\left(\varepsilon\right)$ such
that
\[
\mathrm{d}\left(\tilde{\boldsymbol{x}}_{t},A\right)\le\frac{\varepsilon}{L_{T}},
\]
and, in case of asymptotic stability, such that
\[
\tilde{\boldsymbol{x}}_{t}\to A.
\]
This means, we can represent $\tilde{\boldsymbol{x}}_{t}$ as $\tilde{\boldsymbol{x}}_{t}=\boldsymbol{a}_{t}+\boldsymbol{\chi}_{t}$
with $\boldsymbol{a}_{t}\in A$, $\left|\boldsymbol{\chi}_{t}\right|\le\frac{\varepsilon}{L_{T}}$
and, in case of asymptotic stability, $\left|\boldsymbol{\chi}_{t}\right|\to0$,
for $t\to\infty$. Note that $t\mapsto\boldsymbol{a}_{t}$ has not
to be a solution. Define $\boldsymbol{b}_{t}=T\left[\boldsymbol{a}_{t}\right]\in T\left[A\right]$.
Then,
\begin{eqnarray*}
\left|\tilde{\boldsymbol{y}}_{t+\bar{\eta}}-\boldsymbol{b}_{t}\right| & = & \left|T\circ\tilde{T}\left[\tilde{\boldsymbol{y}}_{t}\right]-T\left[\boldsymbol{a}_{t}\right]\right|\\
 & = & \left|T\left[\boldsymbol{a}_{t}+\boldsymbol{\chi}_{t}\right]-T\left[\boldsymbol{a}_{t}\right]\right|\\
 & \le & L_{T}\left|\boldsymbol{\chi}_{t}\right|\le\varepsilon
\end{eqnarray*}
and $\left|\tilde{\boldsymbol{y}}_{t+\bar{\eta}}-\boldsymbol{b}_{t}\right|\to0$,
if $A$ is asymptotically stable. This completes the proof.
\end{proof}

\section{Reduction of delay-parameters\label{sec:Reduction-of-Delay-Parameters}}

The Theorems \ref{thm:correspondence} and \ref{thm:stability} show
that CT-equivalence is indeed a very strong equivalence. For virtually all
cases one is best advised to study the system which possesses the most
convenient distribution of delays within the concerned equivalence
class. However, it is difficult give a general identification of 
the appropriate distribution, since it strongly depends on the problem
at hand. We can acknowledge two possible guiding principles, which
correspond to the transformations of delays in a ring that were
mentioned in the introduction [cf. (\ref{eq:Baldi-delays}) and (\ref{eq:Perlikowski-delays})].
Firstly, a homogenization of delays can sometimes lead to a higher 
degree of symmetry and thereby allow for simplifications. 
In other cases a "concentration" of the delays on selected links may be useful.

In Sec.~\ref{sub:Finding-a-SpanningTree} we show,
that it is always possible to find timeshifts $\eta_{j}$, $j=1,...,N$,
such that the number of different delays in system (\ref{eq:1}) reduces
to the cycle space dimension $C=L-N+1$ of the network, where $L=\#{\cal E}$
is the number of links and $N=\#{\cal N}$ is the number of nodes in the network.
Effectively, this means that no more than $C$ delay-parameters have
to be taken into account during investigation {[}see Fig.~\ref{fig:CTT-Example}{]}.
For a connected network this number cannot be reduced further as we 
show in Sec.~\ref{sub:Genericity}.

\subsection{Construction of an instantaneous spanning tree\label{sub:Finding-a-SpanningTree}}

In this section we construct a set of links, called a "spanning tree" (see Def.~\ref{def:spanning-tree}),
on which all connection delays can be eliminated by componentwise timeshifts. Let us firstly 
introduce the necessary notions:
\begin{definition}
A semicycle $c=\left(\ell_{1},...,\ell_{k}\right)$ is a closed path
in the undirected graph which is obtained by dropping the orientation
from all links from the multigraph $({\cal N},{\cal E})$. 
\end{definition}

\begin{definition}\label{def:spanning-tree}
A spanning tree of $({\cal N}, {\cal E})$ is a set of links $S\subseteq {\cal E}$
which contains no semicycles but all nodes, that is $\{s(\ell),t(\ell)\}_{\ell\in S}={\cal N}$.
\end{definition}

\begin{figure}%
\begin{minipage}[u]{0.12\linewidth}
(a) \\initial stage:\\
\includegraphics[scale=0.25]{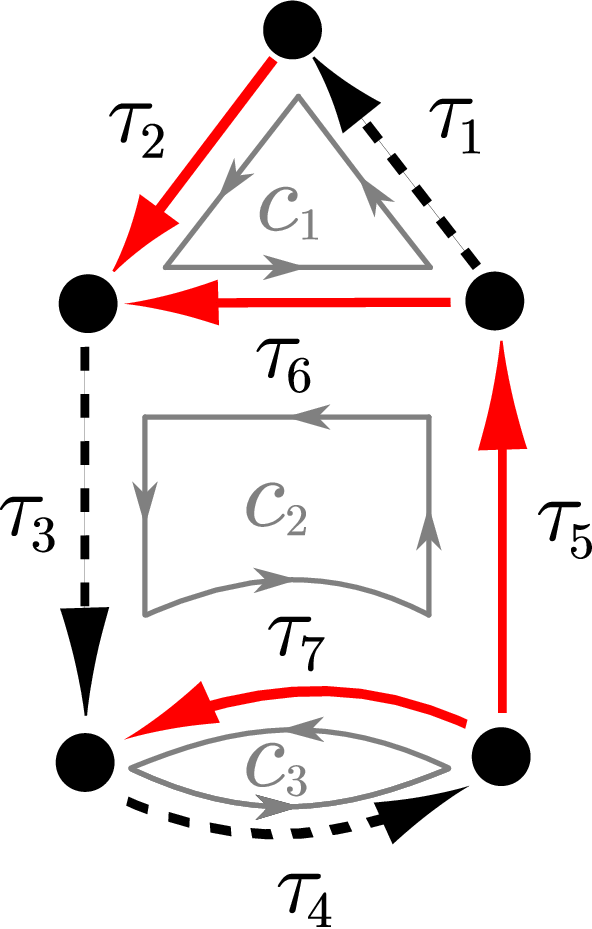} 
\end{minipage}
\begin{minipage}[l]{0.62\linewidth}

\raisebox{-3mm}[0mm][0mm]{(b)}
\begin{tabular}[t]{|r|c|c|c|c|c|c|c|}
\hline 
stage & $\tau_{1}$ & $\tau_{2}$ & $\tau_{3}$ & $\tau_{4}$ & $\tau_{5}$ & $\tau_{6}$ & $\tau_{7}$\tabularnewline
\cline{2-8} 
\raisebox{-2mm}[0mm][0mm]{(i) \raisebox{-1mm}[0mm][0mm]{\rotatebox{90}{$\curvearrowleft$}}} & 4 & 3 & 10 & 8 & 9 & 2 & 5\tabularnewline
\cline{2-8} 
\raisebox{-2mm}[0mm][0mm]{(ii) \raisebox{-1mm}[0mm][0mm]{\rotatebox{90}{$\curvearrowleft$}}} & 4 & \textcolor{red}{{[}3{]}} & 10 & 8 & \textcolor{red}{{[}9{]}} & \textcolor{red}{{[}2{]}} & \textcolor{red}{{[}5{]}}\tabularnewline
\cline{2-8} 
\raisebox{-2mm}[0mm][0mm]{(i) \raisebox{-1mm}[0mm][0mm]{\rotatebox{90}{$\curvearrowleft$}}} & 2 & \textcolor{red}{{[}3{]}} & 12 & 8 & \textcolor{red}{{[}9{]}} & \textcolor{red}{{[}0{]}} & \textcolor{red}{{[}5{]}}\tabularnewline
\cline{2-8} 
\raisebox{-2mm}[0mm][0mm]{(ii) \raisebox{-1mm}[0mm][0mm]{\rotatebox{90}{$\curvearrowleft$}}} & \textcolor{red}{{[}2{]}} & 3 & 12 & 8 & \textcolor{red}{{[}9{]}} & \textcolor{red}{{[}0{]}} & \textcolor{red}{{[}5{]}}\tabularnewline
\cline{2-8} 
\raisebox{-2mm}[0mm][0mm]{(i) \raisebox{-1mm}[0mm][0mm]{\rotatebox{90}{$\curvearrowleft$}}} & \textcolor{red}{{[}0{]}} & 5 & 12 & 8 & \textcolor{red}{{[}9{]}} & \textcolor{red}{{[}0{]}} & \textcolor{red}{{[}5{]}}\tabularnewline
\cline{2-8} 
\raisebox{-2mm}[0mm][0mm]{(ii) \raisebox{-1mm}[0mm][0mm]{\rotatebox{90}{$\curvearrowleft$}}} & \textcolor{red}{{[}0{]}} & 5 & 12 & 8 & \textcolor{red}{{[}9{]}} & \textcolor{red}{{[}0{]}} & \textcolor{red}{{[}5{]}}\tabularnewline
\cline{2-8} 
\raisebox{-2mm}[0mm][0mm]{(i) \raisebox{-1mm}[0mm][0mm]{\rotatebox{90}{$\curvearrowleft$}}} & \textcolor{red}{{[}0{]}} & 5 & 7 & 13 & \textcolor{red}{{[}9{]}} & \textcolor{red}{{[}0{]}} & \textcolor{red}{{[}0{]}}\tabularnewline
\cline{2-8} 
\raisebox{-2mm}[0mm][0mm]{(ii) \raisebox{-1mm}[0mm][0mm]{\rotatebox{90}{$\curvearrowleft$}}} & \textcolor{red}{{[}0{]}} & 5 & \textcolor{red}{{[}7{]}} & 13 & 9 & \textcolor{red}{{[}0{]}} & \textcolor{red}{{[}0{]}}\tabularnewline
\cline{2-8} 
 & \textcolor{red}{{[}0{]}} & 5 & \textcolor{red}{{[}0{]}} & 13 & 16 & \textcolor{red}{{[}0{]}} & \textcolor{red}{{[}0{]}}\tabularnewline
\hline 
\end{tabular}
\end{minipage}
\begin{minipage}[l]{0.2\linewidth}
(c) \\final stage:\\
\includegraphics[scale=0.25]{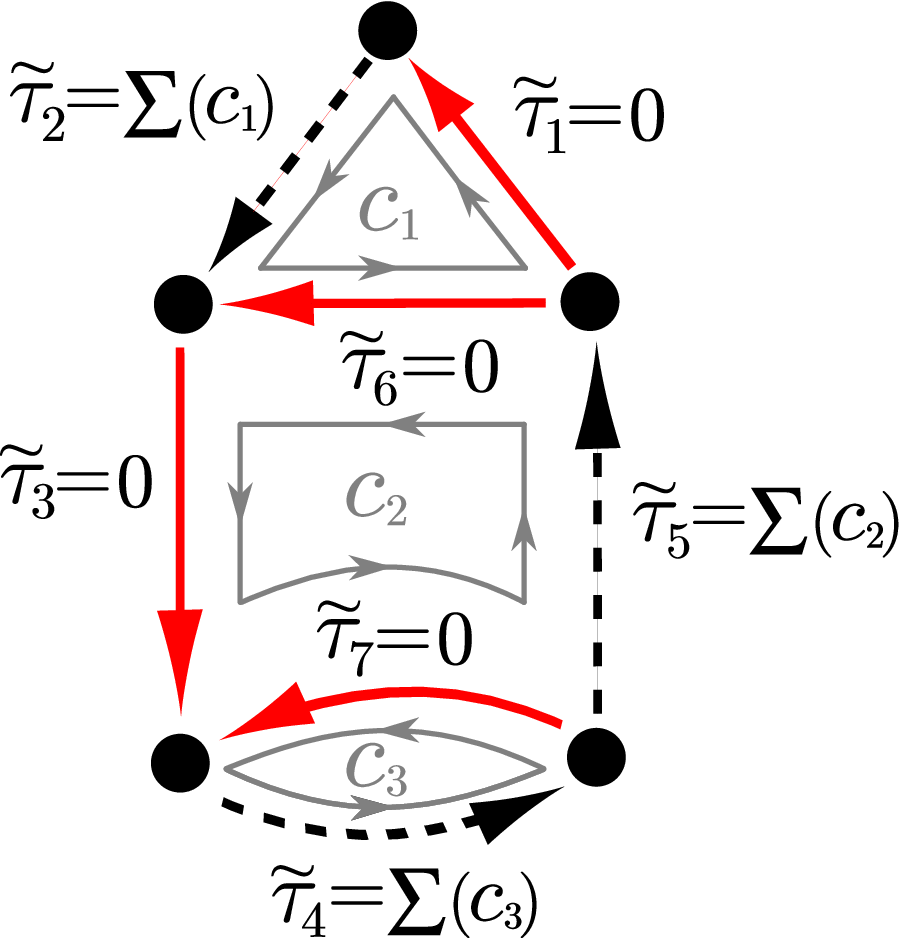}
\end{minipage}
\hspace{1cm}

\par

\caption{\label{fig:spanning-tree-illustration} 
Illustration of the delay reduction in a network of $N=7$ 
delay-coupled systems with cycle space dimension $C=3$ 
[cf. proof of Theorem.~\ref{thm:spanning-tree}]. 
In (a) the initially selected spanning tree with original delays
$\tau_j$, $j=1,...,7$, is indicated by red, solid links,
Three semicycles $c_1$, $c_2$, and $c_3$ are indicated 
by the grey curves following the contained links;
in (c) the final spanning tree with transformed delays $\tilde{\tau}_j$ 
is indicated. The table (b) shows the steps taken by the 
reduction algorithm.
Bracketed, red values indicate the delay times on the
currently selected spanning tree at each step.}
\end{figure}

A spanning tree can be thought of as a "skeleton" of the graph. 
Following its links one can visit each node in the graph exactly once.
For instance, the solid, red links in Figs.~\ref{fig:spanning-tree-illustration}(c) 
form spanning trees. A spanning tree can also be characterized as a maximal cycle-free 
set of links. Consequently, adding another link which is not contained in the spanning 
tree creates a semicycle. For instance, if in Figs.~\ref{fig:spanning-tree-illustration}(a) 
the link with delay $\tau_1$ is added to the spanning tree the semicycle $c_1$ is created.
A spanning tree of a connected network necessarily contains
$N-1$ links, therefore the cycle space dimension $C$ coincides with the number
of links not contained in a spanning tree.
In the context of the CTT, an important quantity for a semicycle is its delay sum.

\begin{definition}
 The delay sum of a semicycle $c=\left(\ell_{1},...,\ell_{k}\right)$ 
with respect to a delay distribution and an orienting link $\ell_1$ is 
\[
\Sigma\left(c\right):=\sum_{j=1}^{k}\sigma_{j}\tau\left(\ell_{j}\right),
\]
where $\sigma_{j}\in\left\{ \pm1\right\} $ indicates whether the
link $\ell_{j}$ points in the same direction as $\ell_{1}$ ($\sigma_{j}=1$)
or not ($\sigma_{j}=-1$). 
The roundtrip of $c$ is the modulus of its delay sum
\begin{equation}
\mathrm{rt}\left(c\right):=\left|\Sigma\left(c\right)\right|.\label{eq:def-roundtrip}
\end{equation}

\end{definition}

The following Lemma gives a characterization of timeshift-transformed delay distributions and describes their relation to the underlying graph structure. More specifically, it shows that the delay sums of semicycles are invariant under timeshifts.
\begin{lemma}
\label{lem:persistence-of-roundtrip} Let $\tau,\tilde{\tau}:{\cal E}\to\left[0,\infty\right)$
be two delay distributions in a connected network. Then, the following
statements are equivalent.

\emph{(i)} There exist $\eta_{j}\ge0$, $j=1,...,N$, such that
$\tilde{\tau}\left(\ell\right)=\tau\left(\ell\right)-\eta_{t\left(\ell\right)}+\eta_{s\left(\ell\right)}.$

\emph{(ii)} For all semicycles $c$ of the network holds:
\[
\Sigma_{\tau}\left(c\right)=\Sigma_{\tilde{\tau}}\left(c\right).
\]
\end{lemma}
\begin{proof}
$\text{(i)}\Rightarrow\text{(ii)}$ Indeed, for any cycle $c=\left(\ell_{1},...,\ell_{k}\right)$
we have
\begin{align*}
\Sigma_{\tilde{\tau}}\left(c\right) &= \sum_{j=1}^{k}\sigma_{j}\tilde{\tau}\left(\ell_{j}\right)
 = \sum_{j=1}^{k}\sigma_{j}\left(\tau\left(\ell_{j}\right)-\eta_{t\left(\ell_{j}\right)}+\eta_{s\left(\ell_{j}\right)}\right)\\
 &= \Sigma_{\tau}\left(c\right)+\sum_{j=1}^{k}\sigma_{j}\left(\eta_{s\left(\ell_{j}\right)}-\eta_{t\left(\ell_{j}\right)}\right)
 = \Sigma_{\tau}\left(c\right).
\end{align*}

$\text{(ii)}\Rightarrow\text{(i)}$ Select an arbitrary spanning tree $S=\left(\ell_{1},...,\ell_{N-1}\right)$.
Select $\xi=\left(\xi_{1},...,\xi_{N}\right)$ and $\chi=\left(\chi_{1},...,\chi_{N}\right)$
such that $\tau\left(\ell_{j}\right)=\xi_{t\left(\ell_{j}\right)}-\xi_{s\left(\ell_{j}\right)}$
and $\tilde{\tau}\left(\ell_{j}\right)=\chi_{t\left(\ell_{j}\right)}-\chi_{s\left(\ell_{j}\right)}$
for $j=1,...,N-1$. (Note that both defining sets of equations are
inhomogeneous linear systems of type $\mathbb{R}^{N\times\left(N-1\right)}$
and of rank $N-1$. Therefore they possess solutions.) The shifts
$\xi$ and $\chi$ define distributions (with possibly negative values)
\begin{eqnarray*}
\hat{\tau}_{1}\left(\ell\right) & = & \tau\left(\ell\right)-\xi_{t\left(\ell\right)}+\xi_{s\left(\ell\right)},\\
\hat{\tau}_{2}\left(\ell\right) & = & \tilde{\tau}\left(\ell\right)-\chi_{t\left(\ell\right)}+\chi_{s\left(\ell\right)}.
\end{eqnarray*}
Both, $\hat{\tau}_{1}$ and $\hat{\tau}_{2}$, are instantaneous along
$S$, i.e., $\hat{\tau}_{1}\left(\ell\right)=\hat{\tau}_{2}\left(\ell\right)=0$
for $\ell\in S$. Each link $\ell$ which is not in $S$ 
corresponds to a unique fundamental semicycle $c=c(S,\ell)$ which is created 
by adding $\ell$ to $S$. It is the only
link in $c$ which may hold a non-zero value $\hat{\tau}_{1}\left(\ell\right)$
or $\hat{\tau}_{2}\left(\ell\right)$, respectively. By the part $\text{(i)}\Rightarrow\text{(ii)}$
we have
\[
\hat{\tau}_{1}\left(\ell\right)=\Sigma_{\tau}\left(c\right)=\Sigma_{\tilde{\tau}}\left(c\right)=\hat{\tau}_{2}\left(\ell\right).
\]
Hence $\hat{\tau}_{1}=\hat{\tau}_{2}$. With $\eta_{j}:=\xi_{j}-\chi_{j}$
this gives
\begin{eqnarray*}
\tilde{\tau}\left(\ell\right) & = & \tau\left(\ell\right)-\eta_{t\left(\ell\right)}+\eta_{s\left(\ell\right)}.
\end{eqnarray*}

Note that $\eta_{j}$ can always be chosen to be non-negative since
a simultaneous shift $\eta_{j}\mapsto\eta_{j}+\bar{\eta}$ by some
amount $\bar{\eta}\in\mathbb{R}$ does not change the resulting transformed
distribution $\tilde{\tau}$.
\end{proof}

Often, this lemma makes it straightforward to determine possible
transformations, since it circumvents the explicit determination of timeshifts.
For instance, in the case of the two coupled systems depicted
in Fig.~\ref{fig:CTT-Example}, it follows immediately that $\tau=(\tau_{1}+\tau_{2})/2$.

Now we can state the main result of this section, which is the construction
of a transformation (\ref{eq:new-variables-y}) such that the number
of delays is minimized.

\begin{theorem}\label{thm:spanning-tree}
For every connected network with dynamics given by Eq.~(\ref{eq:1}),
there exists a spanning tree $S$ and timeshifts $\eta_{j}$ such
that in the transformed system (\ref{eq:1-transformed}) all links
$\ell\in S$ are instantaneous, 
\[
\tilde{\tau}\left(\ell\right)=0,\text{ \emph{for} }\ell\in S,
\]
and for each link $\ell$ outside the spanning tree the delay $\tilde{\tau}\left(\ell\right)$
equals the roundtrip $rt\left(c\left(\ell\right)\right)$ {[}Eq.~(\ref{eq:def-roundtrip}){]}
along the corresponding fundamental cycle $c\left(\ell\right)$, 
\[
\tilde{\tau}\left(\ell\right)=T\left(c\left(\ell\right)\right)\text{, \emph{for} }\ell\notin S.
\]
\end{theorem}
\begin{proof}
The following algorithm describes the general procedure to find
a reduced set of delays. For a more definite example 
see Fig.~\ref{fig:spanning-tree-illustration} and its description following this proof.
The main idea of the algorithm is as follows. We construct spanning
trees and timeshifts iteratively such that after each step the number
of non-zero delays on the spanning tree has decreased by at least
one. Therefore the algorithm finishes after at most $N-1$ iterations.
Each step consists of two stages.

\emph{Stage} (i). Construction of a spanning tree: \\
Select a spanning tree $S=\left\{ \ell_{1},....,\ell_{N-1}\right\} $
in the following way. First, pick a link $\ell_{1}$ with minimal
delay, i.e. $\tau\left(\ell_{1}\right)=\min_{\ell\in{\cal L}}\left\{ \tau\left(\ell\right)\right\} \in[0,\infty)$.
Proceed picking links $\ell_{2},...,\ell_{j}$ with minimal delays
(under all links except the ones already picked) as long as the set
$S_{j}=\left\{ \ell_{1},...,\ell_{j}\right\} $ contains no semicycles.
If at the $j$-th step the chosen link would create a semicycle when
added to $S_{j-1}$, do not add it but ignore it for the rest of this
stage. Following this procedure, called Kruskal's algorithm \cite{Kruskal1956},
yields a spanning tree $S=\left\{ \ell_{1},...,\ell_{N-1}\right\} $.

\emph{Stage} (ii). Construction of timeshifts: \\
Now consider the link $\ell^{\ast}\in S$ with minimal \emph{positive}
delay (in the very first step of the construction it equals $\ell_{1}$
if $\tau\left(\ell_{1}\right)>0$). The link induces a \emph{fundamental
cut} \cite{Diestel2010}, i.e. it partitions the spanning tree $S$
in two connected components: the source component $V$ of $\ell^{\ast}$
and the target component $W$ of $\ell^{\ast}$, see Fig.~\ref{fig:fundamental-section}.
\begin{figure}%
\centering{}\includegraphics[width=0.65\textwidth]{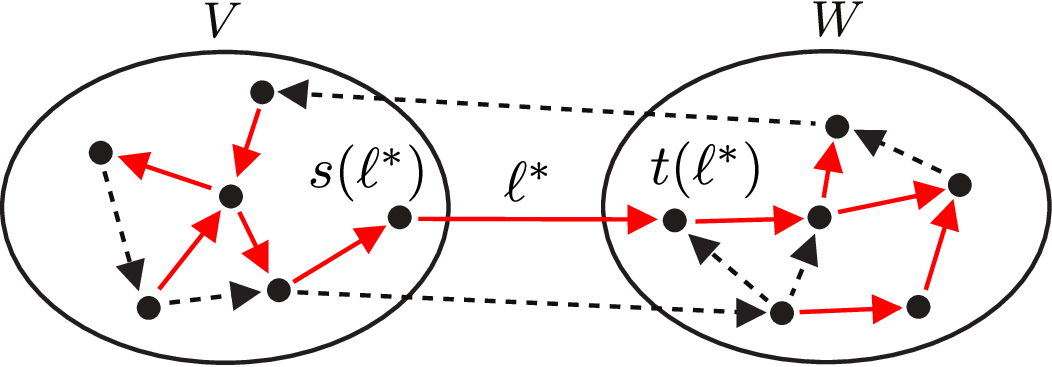}\caption{\label{fig:fundamental-section}
Illustration of the fundamental cut
corresponding to the link $\ell^{\ast}$ in the spanning tree $S$
(solid red links), $s(\ell^\ast)$ and $t(\ell^\ast)$ denote the source
and target of $\ell^\ast$. Links which are not contained in $S$ are indicated
by dashed lines. }
\end{figure}%
Since $S$ is spanning, this is a partition of all nodes.
Now let us define the timeshifts by $\eta_{j}=0$ for $j\in V$ and
$\eta_{j}=\tau\left(\ell^{\ast}\right)$ for all $j\in W$.\\
From the shifts $\eta_{j}$ we obtain the new delay distribution $\tilde{\tau}\left(\ell\right)=\tau\left(\ell\right)+\eta_{s\left(\ell\right)}-\eta_{t\left(\ell\right)}$.
For a link $\ell\ne\ell^{\ast}$ which connects nodes within one set
of the partition, i.e. $s(\ell),\, t(\ell)\in V$ or $s(\ell),\, t(\ell)\in W$,
we have $\eta_{s\left(\ell\right)}=\eta_{t\left(\ell\right)}$. This
implies $\tilde{\tau}\left(\ell\right)=\tau\left(\ell\right)$. In
particular this is the case for all links in $S\setminus\left\{ \ell^{\ast}\right\} $
and therefore $\tau(\ell)=0=\tilde{\tau}(\ell)$ for all instantaneous
links in $S$. For $\ell^{\ast}$ we have $\tilde{\tau}\left(\ell^{\ast}\right)=0$
because $s\left(\ell^{\ast}\right)\in V$ and $t\left(\ell^{\ast}\right)\in W$,
which means $\eta_{s\left(\ell^{\ast}\right)}-\eta_{t\left(\ell^{\ast}\right)}=-\tau\left(\ell^{\ast}\right)$.\\
Hence, the delay distribution $\tilde{\tau}$ reduced the number of
delays on the spanning tree $S$ by one compared to the distribution
$\tau$. We remind that $\tau\left(\ell^{\ast}\right)$ is the smallest
positive delay not only on the spanning tree but over all the links
which connect $V$ and $W$. Therefore, $\tau(\ell)\ge0$ for all
$\ell\in{\cal L}$. Indeed, if there were a link between $V$ and
$W$ with a smaller delay we must have included it in the spanning
tree in step (i)

If a link $\ell$ exists in $S$ with $\tilde{\tau}\left(\ell\right)\ne0$
we repeat stage~(i) starting from an initial set $S_{j}$ which contains
all $j$ instantaneous links from the spanning tree constructed in
the previous step. Then stage~(ii) yields a spanning tree with a
strictly larger number of instantaneous links. In this way we arrive
at an instantaneous spanning tree $S$ in at most $N-1$ iterations.
The statement that $\tilde{\tau}\left(\ell\right)=T\left(c\left(\ell\right)\right)$
for $\ell\notin S$ follows from Lemma~\ref{lem:persistence-of-roundtrip}.
\end{proof}

Fig.~\ref{fig:spanning-tree-illustration} illustrates how the algorithm 
described in the proof of Theorem.~\ref{thm:spanning-tree} determines a delay 
reduction in the depicted network of $N=7$ coupled systems with 
cycle space dimension $C=3$. We assume an initial delay-distribution 
as given in the first row of the table in Fig.~\ref{fig:spanning-tree-illustration}(b). 
Then, the spanning tree in (a), indicated by solid, red links, is selected by the algorithm
as a successive cycle-free collection of links with smallest delays.
Let us denote the links of the example by $\ell_j$ with $\tau_j = \tau(\ell_j)$, $j=1,...,7$.
Note that in (a) $\ell_1$ is not contained in the spanning tree although 
$\tau_1 < \tau_5$ and $\tau_1 < \tau_7$, and $\ell_5$ and $\ell_7$ are both contained.
This is because if the links $\ell_2$ and $\ell_6$ with $\tau_2,\tau_6<\tau_1$ are 
already selected, then an addition of $\ell_1$ would lead to the inclusion 
of the semicycle $c_1$, which is not allowed.
Each row of table (b) corresponds to a step of the algorithm, where 
the bracketed, red values correspond to the delay times of the currently 
selected spanning tree. Two stages are repeated alternatingly:
In stage (i) a spanning tree with minimal delay sum is determined and in stage (ii) the
minimal delay in this spanning tree is eliminated by a componentwise timeshift. 
In (c) the spanning tree is shown which is selected at the final stage. 
It contains only instantaneous links when the algorithm is completed. 
The remaining delays on links not contained in the final spanning tree,
correspond to the delay sums of the corresponding semicycles ($c_1$, $c_2$, $c_3$), 
oriented as indicated in the figure, i.e.,
$\tilde{\tau}_2=\Sigma(c_1)=\tau_1+\tau_2-\tau_6$;
$\tilde{\tau}_5=\Sigma(c_2)=\tau_3+\tau_5+\tau_6-\tau_7$; and
$\tilde{\tau}_4=\Sigma(c_3)=\tau_4+\tau_7$.

\subsection{Genericity of the dimension of the delay parameter space\label{sub:Genericity}}

We call $C$ the \emph{essential number }of delays\emph{ }since it
is the minimal number to which the number of different delays
in a network can be reduced generically. Here {}``generically''
means that the conditions which allow for further reduction of delays
form a null set in the parameter space of delays $\mathbb{R}_{\ge0}^{L}=\{\left(\tau\left(\ell\right)\right)_{\ell\in{\cal E}}\,:\,\tau\left(\ell\right)\ge0\}$
of the original system. The reducibility condition to $m$ different
delays is described by the $L$ linear equations 
\begin{equation}
\eta_{t\left(\ell\right)}-\eta_{s\left(\ell\right)}+\tilde{\tau}\left(\ell\right)=\tau\left(\ell\right),\ \ell\in I_{j},\label{eq:redux-condition}
\end{equation}
with the restriction for $\tilde{\tau}\left(\ell\right)$ to take
one of $m$ different values $\left\{ \theta_{1},...,\theta_{m}\right\}$.
System (\ref{eq:redux-condition}) can be equivalently written in
the vector form $G_{q}\mathbf{v}=\mathbf{\boldsymbol{\tau}}$, where
$\mathbf{v}$ is the $(N+m)$-dimensional vector of unknowns $\mathbf{v}=(\eta_{1},\dots,\eta_{N},\theta_{1},\dots\theta_{m})$
and $\boldsymbol{\tau}$ is the $L$-dimensional vector of delays
$\tau\left(\ell\right)$. For any fixed assignment $\tilde{\tau}\left(\ell\right)=\theta_{q\left(\ell\right)}$,
$\ell\in I_{j}$, with $q:I_{j}\to\left\{ 0,...,m\right\} $, and
$\theta_{0}:=0$, one obtains a different matrix $G_{q}\in\mathbb{R}^{L\times\left(N+m\right)}$.
It can further be shown that the rank of the matrix $G_{q}$ is smaller
than $N-1+m$. This implies for $m<C=L-(N-1)$ that $N-1+m<L$. Hence,
the number of equations in (\ref{eq:redux-condition}) is larger than
the number of unknowns. Such equation cannot be solved generically,
unless the given delays $\tau\left(\ell\right)$ satisfy some special
condition of positive codimension.

\section{Conclusion}

In this article we have studied a componentwise timeshift transformation
(CTT) for a general class of coupled differential equations with constant
coupling delays. We have defined appropriate phase spaces such that
the CTT conveys an equivalence for the semiflows of the original and
the transformed system which is reminiscent to, but weaker than topological
conjugacy. We have shown several dynamical invariants for the equivalent
flows. Firstly, for delay differential equations, the characteristic
exponents of equilibria and the Floquet exponents of periodic orbits
are invariant under the CTT. More generally, for CT-equivalent semidynamical
systems, we have shown that there exists a one-to-one correspondence
between invariant and strongly invariant sets, respectively. As a
main stability result we have shown that corresponding positively
invariant sets have the same kind of stability, and invariant and
strongly invariant corresponding sets possess the same maximal Lyapunov
exponents. To sum up, the observable dynamics of the coupled units
might change its relative timing in the transformed system, but qualitatively
it remains the same. In particular attractors and their stability
are invariant.

We have presented a constructive proof that there is always a CTT
which reduces the delays to at most $C=L-N+1$ different delays, where
$L$ is the number of links in the network and $N$ is the number
of nodes. The number of different delays $C$ coincides with the cycle
space dimension of the underlying graph. Furthermore, we have shown
that the sum of delays along semicycles in the network is invariant
under CTTs. This reveals a strong link between a coupled delay differential
equation and the topology of the underlying graph. Although a CTT
which reduces the number of different delays to a minimum is usually
not unique, we have shown that the minimal number of different delays
itself cannot be reduced in general.

We believe that our results have a relevance for several applied areas
(see also \cite{Luecken2013b}). For example, in theoretical neuroscience
it is an accepted fact that in the case of two mutually delay coupled
units \cite{Panchuk2013} (as in Fig.~\ref{fig:CTT-Example}) or
more generally, in unidirectionally coupled rings \cite{Perlikowski2010},
the delays can always be chosen identical for theoretical investigations.
Moreover, in view of our results, the observations about the role
of the greatest common divisor of loop-lengths in networks of delay
coupled excitable systems with homogeneous delays \cite{Kanter2011,Rosin2013}
can be formulated more generally in terms of the delay sums along
the network's cycles. The CTT is an important tool for investigators
working with delayed dynamical systems since it clarifies one aspect
of the way in which different interaction delays in a coupled system
work together.

Furthermore, the CTT can speed up the numerical simulation of a system
by reducing the number of different delays or by reducing the maximal
delay time. The computational advantages of the transformed system
might be of interest especially in the field of delayed neural networks,
where large scale simulations of networks with many different delay
times are conducted.

\bibliographystyle{plain}

\end{document}